\documentclass[11pt,leqno,a4paper]{article}
\usepackage{amsmath,amssymb,amsthm,amsxtra,latexsym}
\usepackage{upref}

\newtheorem{corollary}{Corollary}[section]
\newtheorem{lemma}{Lemma}[section]
\newtheorem{proposition}{Proposition}[section]
\newtheorem{theorem}{Theorem}[section]
\theoremstyle{definition}
\newtheorem{remark}{Remark}
\newtheorem{definition}{Definition}[section]

\theoremstyle{remark}

\newcommand{\R}{\mathbb{R}}
\newcommand{\C}{\mathbb{C}}
\newcommand{\N}{\mathbb{N}}
\newcommand{\Z}{\mathbb{Z}}
\newcommand{\la}{\langle}
\newcommand{\ra}{\rangle}
\newcommand{\pd}{\partial}

\newcommand{\bM}{\mathbb{M}}

\title{Asymptotic stability of small solitons for 2D Nonlinear
Schr\"{o}dinger equations with potential}
\author{Tetsu Mizumachi
\footnote{This research is supported by Grant-in-Aid for Scientific Research
(No. 17740079).}}
\date{}
\begin{document}
\maketitle
\begin{abstract}
We consider asymptotic stability of a small solitary wave to
supercritical $2$-dimensional nonlinear Schr\"{o}dinger equations
$$ iu_t+\Delta u=Vu\pm |u|^{p-1}u \quad\text{for 
$(x,t)\in\mathbb{R}^2\times\mathbb{R}$,}$$
in the energy class. This problem was studied by
Gustafson-Nakanishi-Tsai \cite{GNT} in the $n$-dimensional case
$(n\ge 3)$ by using the endpoint Strichartz estimate.
Since the endpoint Strichartz estimate fails in 2-dimensional
case, we use a time-global local smoothing estimate of Kato type
to prove the asymptotic stability of a solitary wave.
\end{abstract}
\maketitle

\section{Introduction}
\label{intro}
In this paper, we consider asymptotic stability of solitary wave
solutions to
\begin{equation}
  \label{eq:NLS}
\left\{
  \begin{aligned}
& iu_t+\Delta u=Vu+f(u) \quad\text{for 
$(x,t)\in\mathbb{R}^2\times\mathbb{R}$},
\\
& u(x,0)=u_0(x)\quad\text{for $x\in\mathbb{R}^2$,}
  \end{aligned}\right.
\end{equation}
where $V(x)$ is a real potential, $f(u)=\alpha|u|^{p-1}u$ with $\alpha=\pm1$.

Let 
\begin{align*}
& H(u)=\int_{\R^2}\left(|\nabla u|^2+V(x)|u|^2
+\frac{2\alpha}{p+1}|u|^{p+1}\right)dx,  
\\
& N(u)=\int_{\R^2}|u|^2dx.
\end{align*}
Then a solution to \eqref{eq:NLS} satisfies
\begin{equation}
  \label{eq:cons}
H(u(t))=E(u_0),\quad N(u(t))=N(u_0)  
\end{equation}
during the time interval of existence. Stability of solitary waves was
first studied by Cazenave and Lions \cite{CaLi},
Grillakis-Shatah-Strauss \cite{GSS} and Weinstein \cite{W2} 
(see also Rose-Weinstein \cite{RW}, Oh \cite{Oh} and Shatah-Strauss
\cite{ShSt}). In the case of integrable equations such as cubic NLS
and KdV, the inverse scattering theory tells us that if the initial data
decays rapidly as $x\to\pm\infty$, a solution decomposes into
a sum of solitary waves and a radiation part as $t\to\infty$
(see \cite{Sc}). Soffer and Weinstein \cite{SW1,SW2}
considered 
\begin{equation}
  \label{eq:NLS2}
iu_t+\Delta u=Vu\pm |u|^{p-1}u \quad\text{for $x\in\R^n$ and $t>0$,}  
\end{equation}
where $n\ge 2$ and $1<p<(n+2)/(n-2)$. They proved that if $-\Delta +V$
has exactly one eigenvalue with negative value $E_*$ and initial data
is well localized and close to a nonlinear bound state, a solution tends
to a sum of a nonlinear bound state nearby and a radiation part which
disperses to $0$ as $t\to \infty$ 
(see also \cite{Kir} for 2-dimensional case).
This result was extended by Yau and Tsai \cite{TY,YT1,YT2,YT3} and 
Soffer-Weinstein \cite{SW3} to the case where $-\Delta+V$ have two bound
states (see also \cite{FTY,Ts}).
In the $1$-dimensional case, Buslaev and Perelman \cite{BP1, BP2}
and Buslaev and Sulem \cite{BuS} studied the asymptotic stability of
\eqref{eq:NLS} with $V\equiv0$.
Using the Jost functions, they built a local energy decay estimate of
solutions to the linearized equation and prove asymptotic stability 
of solitary waves for super critical nonlinearities. 
Their results are extended to the higher dimensional case by
Cuccagna \cite{Cu}.  See also Perelman \cite{Per} and 
Rodnianski-Schlag-Soffer \cite{RSS} which study asymptotic stability
of multi-solitons, and Krieger and Schlag \cite{KrSc} which study large time
behavior of soluitons around unstable solitons.
\par
However, all these results assume that initial data is well localized
 so that a solution decays like $t^{-3/2}$.
Martel and Merle \cite{MM}, \cite{MM2} proved the asymptotic stability
of solitary waves to generalized KdV equations using the monotonicity of $L^2$-mass, which is a variant of the local smoothing effect proved by
Kato \cite{K}. They elegantly use the fact that the dispersive remainder
part of a solution $v(t,x)$ satisfies
\begin{equation}
  \label{eq:disp}
\int_0^\infty \|v(t,\cdot)\|_{H^1_{loc}}^2dt<\infty  
\end{equation}
to prove the asymptotic stability of solitary waves in $H^1$
(see also Pego and Weinstein \cite{PW2} for KdV with exponentially localized
initial data and Mizumachi \cite{Mi1} for polynomially localized initial
data). Gustafson-Nakanishi-Tsai \cite{GNT} has proved asymptotic
stability of a small solitary wave of \eqref{eq:NLS2} in the energy
class with $n\ge3$. Their idea is to use the endpoint Strichartz
estimate instead of \eqref{eq:disp}, which tells us that
$\|v\|_{L^2_tW^{1,6}_x}$ remains small globally in time
for super critical nonlinearities.
However, dispersive wave decays more slowly in the lower dimensional
case and the endpoint Strichartz estimate does not hold
in the $2$-dimensional case. 
Recently, Mizumachi \cite{Mi} has proved the asymptotic stability of
small solitons in 1D case by using dispersive estimates such as
\begin{equation}
\label{eq:tetsu}
 \|\pd_xe^{{it(-\pd_x^2+V)}}P_cf\|_{L_x^\infty L^2_t}
\le C\|f\|_{H^{1/2}}.
\end{equation}
In the present paper, we apply a local
smoothing estimate
\begin{equation}
  \label{eq:disp2}
 \|\la x\ra^{-1-0}e^{it(-\Delta+V)}P_cf\|_{L^2_t(0,\infty;L_x^2(\R^2))}
\le C\|f\|_{L^2(\R^2)},
\end{equation}
to obtain the asymptotic stability of small solitons in the
$2$-dimensional case. 
\par
Local smoothing estimates such as \eqref{eq:disp2} have been studied
by many authors. See, for example,  Ben-Artzi and Klainerman
\cite{Be-Kl}, Constantin and Saut \cite{CS}, Kato and Yajima \cite{KY}
and Kenig-Ponce-Vega \cite{KPV1,KPV2}, Sjolin \cite{Sj}, Ruiz-Vega
\cite{RZ1}, Sugimoto \cite{Sug} and Watanabe \cite{Watanabe}.
Especially, Ben-Artzi and Klainerman \cite{Be-Kl} and
Barcel\'{o}-Ruiz-Vega \cite{BRVe} prove time-global local
smoothing estimates in $n$-dimensional case with $n\ge3$. 
In the $2$-dimensional case, it is well-known that
\begin{equation}
  \label{eq:locsm1}
\|e^{it\Delta} f\|_{L^\infty_x(\R^2;L^2_t(\R))}
\lesssim \|f\|_{L^2(\R^2)},
\end{equation}
follows from a special case of Thomas-Stein theorem (\cite{Ste})
(see, e.g., Planchon \cite{Pl}).
However, to the best of our knowledge, 
there seems to be a lack of literature
in the $2$-dimensional case with $V\not\equiv0$.
 Another purpose of the present is to fill the gap.
\par
Our strategy to prove \eqref{eq:disp2} is to apply Plancherel's theorem
to the inversion of Laplace formula. The key is to prove
\begin{equation}
  \label{eq:res-est}
\|\la x\ra^{-1-0}R(\lambda\pm i0
)f\|_{L^2_\lambda(0,\infty;L^2_x)}
\le C\|f\|_{L^2} \quad\text{for every $f\in L^2(\R^2)$.}  
\end{equation}
To obtain \eqref{eq:res-est}, we prove  that the free resolvent
operator $R_0(\lambda)=(-\Delta-\lambda)^{-1}$ satisfies
\begin{equation}
  \label{eq:frres}
\sup_x\|R_0(\lambda\pm i0)f\|_{L^2_\lambda(0,\infty)}\le C\|f\|_{L^2}
\quad\text{for every $f\in L^2(\R^2)$,}
\end{equation}
and apply a resolvent expansion obtained by Jensen and Nenciu \cite{JeN}
as well as Schlag \cite{Sch}.
We remark that, roughly speaking, Eq. \eqref{eq:frres} can be
translated into \eqref{eq:locsm1} by using the Fourier transform with
respect to $\lambda$.
\par
Our plan of the present paper is as follows.
In Section 2, we state our main result and linear dispersive estimates
that will be used later. In Section 3, we prove our main result
assuming the linear estimates introduced in Section 2.
 In Section 4, we prove \eqref{eq:res-est} and obtain
\eqref{eq:disp2}. To prove \eqref{eq:res-est}, we use an argument
of the resolvent expansion as well as \eqref{eq:frres} which 
follows from $L^2(0,\infty;\sqrt{x}dx)$-boundedness of the Hankel
transform and the $\mathcal{Y}_0$-transform (see Rooney \cite{R3}).

\par
 Finally, we introduce several notations.
Let
\begin{gather*}
\|f\|_{L_t^qL_x^p}=\biggl(\int_\R\bigl(\int_{\R^2}|f(t,x)|^pdx\bigr)^{q/p}dt
\biggr)^{1/q},
\\
\|f\|_{L_x^sL_t^r}=\biggl(\int_{\R^2}\bigl(\int_\R |f(t,x)|^rdt\bigr)^{s/r}dx
\biggr)^{1/s}.
\end{gather*}
We denote by $L^{2,s}$ and $H^{m,s}$ Hilbert spaces whose norms
are defined by  $$\|u\|_{L^{2,s}}=\|\la x\ra^su\|_{L^2(\R^2)}
\quad\text{and}\quad\|u\|_{H^{m,s}}=\|\la x\ra^su\|_{H^m(\R^2)},$$
 where $m\in \N$, $s\in \R$ and $\la x\ra=(1+|x|^2)^{1/2}$.
Let
$$\la f_1,f_2\ra_x=\int_{\R^2} f_1(x)f_2(x)dx,\quad
\la g_1,g_2\ra_{t,x}=\int_{-\infty}^\infty\int_{\R^2}
g_1(t,x)g_2(t,x)dxdt.$$
We set $L^2_{rad}=\{f\,|\, f\in L^2(0,\infty;rdr)\}$ whose norm is
defined by
$$\|f\|_{L^2_{rad}}=(\int_0^\infty |f(r)|^2rdr)^{1/2}.$$
For any Banach spaces $X$, $Y$, we denote by $B(X,Y)$ the space of
bounded linear operators from $X$ to $Y$. We abbreviate
$B(X,X)$ as $B(X)$.

We define the Fourier and transform of $f(x)$ as
$$
\mathcal{F}_xf(\xi)=
(2\pi)^{-n/2}\int_{\R^n} f(x)e^{-ix\xi}dx,$$
and the inverse Fourier transform of $g(\xi)$ as
$$
\mathcal{F}^{-1}_\xi g(x)=
(2\pi)^{-n/2}\int_{\R^n} g(\xi)e^{ix\xi}d\xi.$$
We denote by $\mathcal{S}_t\otimes\mathcal{S}_x(\R^2)$ a set of
functions $f(t,x)=\sum_{i=1}^N f_i(t)g_i(x)$ with $f_i\in \mathcal{S}(\R)$,
$g_i\in \mathcal{S}(\R^2)$ ($1\le i\le N$).

\section{The Main result and Preliminaries}
In the present paper, we assume that the linear potential $V(x)$ is a
$C^1$-function on $\R^2$ satisfying the following.
\begin{itemize}
\item[(H1)] There exists a $\sigma>3$ such that
$\sup_{x\in\R^2}\left(\la x\ra^\sigma|V(x)|
+|\nabla V(x)|\right)<\infty$.
\item [(H2)] $L=-\Delta+V$ has exactly one negative eigenvalue $E_*$.
\item [(H3)] $0$ is neither a resonance nor an eigenvalue of $L$
(see Definition \ref{def:NR} in Section \ref{sec:4}).
\end{itemize}
From (H1)--(H3), it follows that the spectrum of $L$ consists of the
continuous spectrum $\sigma_{c}(L)=[0,\infty)$ and an discrete
eigenvalue $E_*$, and that  $\lambda=E_*$ is a simple eigenvalue of
$L$ (see \cite{RS}). Let $\phi_*$ be a normalized eigenfunction of $L$
(satisfying $\|\phi_*\|_{L^2}=1$) belonging to $E_*$, and let $P_{d}$
 and $P_c$ be spectral projections of $L$ defined by
$$P_{d}u=\left\la u, \phi_*\right\ra\phi_*, \quad
P_cu=(I-P_{d})u.$$
\par
Suppose that $E\in\R$ and $e^{-iEt}\phi_E(x)$ be a solitary wave
 solution of
\eqref{eq:NLS}.
Then $\phi_E(x)$ is a solution to
\begin{equation}
\label{eq:B}
\left\{\begin{aligned}
&\Delta\phi_E+E\phi_E=V\phi_E+\alpha|\phi_E|^{p-1}\phi_E 
\quad\text{for $x\in \R^2$},\\
& \lim_{|x|\to\infty}\phi_E(x)=0.
\end{aligned}\right.
\end{equation}

Using the bifurcation theory, we have the following.
\begin{proposition}
  \label{prop:1} Assume (H1)--(H3). Let $\delta$ be a small positive
number. Suppose  that $E\in(E_*,E_*+\delta)$ and $\alpha=1$ or
$E\in (E_*-\delta,E_*)$ and $\alpha=-1$. Then,
there exists a positive solution $\phi_E$ to \eqref{eq:B} such that
for every $k\in\N$,
\begin{enumerate}
\item $\phi_E\in H^{1,k}$,
\item The function $E\mapsto \phi_E$ is $C^1$ in
$H^{1,k}$ for every $k\in\N$, and as $E\to E_*$,
$$
\phi_E=|E-E_*|^{1/(p-1)}\left(\|\phi_*\|_{L^{p+1}}^{-(p+1)/(p-1)}
\phi_*+O(E-E_*)\right)\quad\text{in $H^{1,k}$}.$$
\end{enumerate}
\end{proposition}
Proposition \ref{prop:1} follows from a rather standard argument.
See for example \cite{Ni} and \cite[pp.123--124]{SW1}.  
\begin{remark}
\label{rem:1}
Let
$\phi_{1,E}=\|\phi_E\|_{L^2}^{-1}\phi_E$ and
$\phi_{2,E}=\|\pd_E\phi_E\|_{L^2}^{-1}\pd_E\phi_E$.
By Proposition \ref{prop:1},
\begin{equation*}
\|\phi_{1,E}-\phi_*\|_{H^{1,k}(\R)}
+\|\phi_{2,E}-\phi_*\|_{H^{1,k}(\R)}\lesssim |E-E_*|.
\end{equation*}
\end{remark}
\par

Now, we introduce our main result.
\begin{theorem}
\label{thm1} Assume (H1)--(H3).
Let $p\ge3$ and let $\varepsilon_0$ be a sufficiently small positive
number.  Suppose $\|u_0\|_{H^1}<\varepsilon_0$. Then there exist an
$E_+<0$, a $C^1$ real-valued function
$\theta(t)$ and $v_+\in P_cH^1(\R^2)$ such that

$$\lim_{t\to\infty}\|u(t)-e^{i\theta(t)}\phi_{E_+}
-e^{-itL}v_+\|_{H^1(\R^2)}=0.$$
\end{theorem}
\begin{remark}  
Let us decompose a solution to \eqref{eq:NLS} into a solitary wave part and
a radiation part:
\begin{equation}
  \label{eq:1.3}
  u(t,x)=e^{-i\theta(t)}(\phi_{E(t)}(x)+v(t)).
\end{equation}
If we take initial data in the energy class, the dispersive part of
the solutions decays more slowly than they does for well localized
initial data. So, being different from Soffer-Weinstein \cite{SW1,SW2}
or Buslaev-Perelman \cite{BP1}, we cannot expect that
$\int_t^\infty\dot{E}(s)ds$ is integrable.
Thus in general, we need dispersive estimates for a time-dependent
linearized equations to prove asymptotic stability of solitary waves
in $H^1$. To avoid this difficulty, we assume the smallness of
solitary waves so that a generalized kernel of the linearized operator
is well approximated by a $1$-dimensional subspace
$\{\beta \phi_*\,|\,\beta\in\C\}.$
\end{remark}
Substituting \eqref{eq:1.3} into \eqref{eq:NLS}, we obtain
\begin{equation}
  \label{eq:v}
  iv_t=Lv+g_1+g_2+g_3+g_4,
\end{equation}
where 
\begin{align*}
& g_1(t)=-\dot{\theta}(t)v(t), \quad
 g_2(t)=(E(t)-\dot{\theta}(t))\phi_{E(t)}-i\dot{E}(t)\pd_E\phi_{E(t)},
\\ 
& g_3(t)=f(\phi_{E(t)}+v(t))-f(\phi_{E(t)})-\pd_\varepsilon
 f(\phi_{E(t)}+\varepsilon v(t))|_{\varepsilon=0},
\\
& g_4(t)=\pd_\varepsilon f(\phi_{E(t)}+\varepsilon v(t))|_{\varepsilon=0}
=\alpha\phi_{E(t)}^{p-1}\left(\frac{p+1}{2}v(t)
+\frac{p-1}{2}\overline{v(t)}
\right).
\end{align*}

\label{sec:2}
To fix the decomposition \eqref{eq:1.3}, we assume
\begin{equation}
  \label{eq:2.1}
\left\la \Re v(t), \phi_{E(t)}\right\ra
=\left\la \Im v(t),\pd_E\phi_{E(t)}\right\ra=0.
\end{equation}
By Proposition \ref{prop:1}, we have
\begin{equation}
  \label{eq:3.3ap}
|E(0)-E_*|^{1/(p-1)}+\|v(0)\|_{H^1} \lesssim \|u_0\|_{H^1}.  
\end{equation}
Since $u\in C(\R;H^1(\R^2))$, it follows from
the implicit function theorem that there exist a $T>0$ and $E$,
$\theta\in C^1([-T,T])$ such that \eqref{eq:2.1} holds for $t\in[-T,T]$.
See, for example, \cite{GNT} for the proof.

Differentiating \eqref{eq:2.1} with respect to $t$ and substituting
\eqref{eq:v} into the resulting equation, we obtain
\begin{equation}
  \label{eq:2.2}
\mathcal{A}(t)
\begin{pmatrix}  \dot{E}(t)\\ \dot{\theta}(t)-E(t)\end{pmatrix}
=\begin{pmatrix}
\la \Im g_3(t),\phi_{E(t)}\ra \\  \la \Re g_3(t),\pd_E\phi_{E(t)}\ra
\end{pmatrix},\end{equation}
where
\begin{align*}
& \mathcal{A}(t)=\\
& \begin{pmatrix}
 \la\pd_E\phi_{E(t)},\phi_{E(t)}\ra-\la \Re v(t),\pd_E\phi_{E(t)}\ra
& \la \Im v(t),\phi_{E(t)}\ra\\
 \la \Im v(t),\pd_E^2\phi_{E(t)}\ra
& \la\pd_E\phi_{E(t)},\phi_{E(t)}\ra+\la \Re v(t),\pd_E\phi_{E(t)}\ra
\end{pmatrix}.  
\end{align*}

To prove our main result, we will use the Strichartz estimate and
the local smoothing effect of Kato type that is global in time.
The Strichartz estimate follows from $L^\infty$-$L^1$ estimate for
$2$-dimensional Schr\"{o}dinger equations with linear potential 
obtained by Schlag \cite{Sch}. See, for example, \cite{Kl-Tao}.
We say that $(q,r)$ is \textsl{admissible} if $q$ and $r$ satisfy
$2<q\le \infty$, $2\le r<\infty$ and $1/q+1/r=1/2$.
For any $p\in[1,\infty]$, we denote by $p'$ a H\"{o}lder conjugate
exponent of $p$.
\begin{lemma}[Strichartz estimate]
  \label{lem:2.1} Assume (H1)--(H3).
  \begin{itemize}
  \item [(a)]
Suppose that $(q,r)$ is admissible. Then there exists a positive number
 $C$ such that for every $f\in L^2(\R)$,
$$\|e^{-itL}P_cf\|_{L_t^qL_x^r}\le C\|f\|_{L^2}.$$
Furthermore, it holds that
\begin{gather*}
\left\|\int_\R e^{isL}P_cg(s,\cdot)ds\right\|_{L^2_x}
\le C\|g\|_{L_t^{q'}L_x^{r'}}.
\end{gather*}
\item[(b)]
Suppose that $(q_1,r_1)$ and $(q_2,r_2)$ are admissible.
Then there exists a positive number $C$
such  that for every $g(t,x)\in \mathcal{S}(\R\times\R^2)$,
\begin{gather*}
\left\|\int_{0}^te^{-i(t-s)L}P_cg(s,\cdot)ds
\right\|_{L_t^{q_1}L_x^{r_1}}
\le C\|g\|_{L_t^{q_2'}L_x^{r_2'}}.
\end{gather*}
  \end{itemize}
\end{lemma}
Since Lemma \ref{lem:2.1} (a) does not hold with $q=2$,
we use the following local estimate
to show that $dE/dt$ is integrable with respect to $t$.
\begin{lemma}
  \label{lem:2.2}
Assume (H1)--(H3).
Let $s>1$. Then there exists a positive constant $C$ such that
\begin{equation}
  \label{eq:2.2.1}
 \|e^{-itL}P_cf\|_{L^2_tL_x^{2,-s}}\le C\|f\|_{L^2}, 
\end{equation}
for every $f\in \mathcal{S}(\R^2)$ and that
\begin{equation}
\label{eq:2.2.2}
\left\|\int_{\R} e^{isL}P_cg(s,\cdot)ds\right\|_{L^2_x}
\le C\|g\|_{L_t^2L_x^{2,s}},
\end{equation}
for every $g(t,x)\in \mathcal{S}(\R\times\R^2)$.
\end{lemma}
\begin{lemma}
  \label{lem:2.3}
  Let $s>1$. Then there exists a positive constant $C$ such that
\begin{equation}
\label{eq:2.4.1}
\left\|\int_0^t e^{-i(t-s)L}P_cg(s,\cdot)ds
\right\|_{L_t^2L^{2,-s}_x}
\le C\|g\|_{L_t^2L_x^{2,s}}.
\end{equation}
for  every $g(t,x)\in \mathcal{S}(\R^2)$ and $t\in\R$.
\end{lemma}
Since the linear term $g_4$ in \eqref{eq:v} may not belong to
$L^{q'}_tL_x^{r'}$ for admissible $(q,r)$ (because $(q_2,r_2)=
(2,\infty)$ is not admissible), we cannot apply Lemma \ref{lem:2.1} (b)
to $g_4$. Instead, we will use the following to deal with $g_4$. 
\begin{corollary}
  \label{cor:2.4}
Let $(q,r)$ be admissible and let $s>1$. Then there exists a positive number
$C$ such that
\begin{equation}
\label{eq:2.3.1}
\left\|\int_{\R} e^{-i(t-s)L}P_cg(s,\cdot)ds\right\|_{L_t^qL_x^r}
\le C\|g\|_{L_t^2L_x^{2,s}}
\end{equation}
for every
$g(t,x)\in \mathcal{S}(\R\times \R^2)$  and $t\in\R$.
\end{corollary}
Using a lemma by Christ and Kiselev \cite{ChK}, we see that
Corollary \ref{cor:2.4} immedaiately follows from Lemmas \ref{lem:2.1}
and \ref{lem:2.2} (see \cite{SmS}).

The proof of Lemmas \ref{lem:2.2}, \ref{lem:2.3} and
Corollary \ref{cor:2.4} will be given in Section \ref{sec:4}.
\section{Proof of Theorem \ref{thm1}}
\label{sec:3}
In this section, we will prove Theorem \ref{thm1}.
To eliminate $g_1$ in \eqref{eq:v}, we put
\begin{equation}
  \label{eq:trans}
w(t)=e^{-i\theta(t)}v(t).
\end{equation}
Then  \eqref{eq:v} is translated into the integral equation
\begin{equation}
  \label{eq:v'}
w(t)=e^{-itL}w(0)-i\sum_{2\le j\le 4}
\int_0^t e^{-i(t-s)L} e^{-i\theta(s)}g_j(s)ds.
\end{equation}
All nonlinear terms in \eqref{eq:v'} can be estimated in terms of
the following. 
\begin{align*}
& \bM_1(T)=\sup_{0\le t \le T}|E(t)-E_*|,\quad
\bM_2(T)=\|\la x\ra^{-s}P_cw\|_{L_t^2(0,T;H_x^1)},
\\ &
\bM_3(T)=\|\la x\ra^{-s}P_{d}w\|_{L_t^2(0,T;H_x^1)},
\\ & 
\bM_4(T)=\sup_{0\le t\le T}\|P_cw(t)\|_{H^1}
+\|P_cw\|_{L_t^q(0,T;W_x^{1,2p})},
\\ &
\bM_5(T)=\sup_{0\le t\le T}\|P_{d}w(t)\|_{H^1}
+\|P_{d}w\|_{L_t^q(0,T;W_x^{1,2p})}.
\end{align*}
where $2/q=1-1/p$.
\begin{proof}[Proof of Theorem \ref{thm1}]
By Proposition \ref{prop:1}, Remark \ref{rem:1} and \eqref{eq:2.1},
$$\la \pd_E\phi_E,\phi_E\ra=O(|E-E_*|^{2/(p-1)-1}),\quad
|\la v, \pd_E^i\phi_E\ra|\lesssim |E-E_*|^{p/(p-1)-i}\|v\|_{L^2}.$$
Thus by \eqref{eq:2.2}, we have
\begin{align}
\label{eq:3.1}
&|\dot{\theta}(t)-E(t)|
 \lesssim \|\phi_{2,E(t)}v^2\|_{L^1}+\|\phi_{2,E(t)}f(v)\|_{L^1},
\\
\label{eq:3.2}
& |\dot{E}(t)|\lesssim \|\phi_{1,E(t)}v^2\|_{L^1}
+\|\phi_{1,E(t)}f(v)\|_{L^1}.
\end{align}

Suppose that the decomposition \eqref{eq:1.3} with
\eqref{eq:2.1} persists for $0\le t\le T$ and that
$\bM_i(T)$ $(1\le i \le 5)$ are bounded.
Eqs. \eqref{eq:3.1}--\eqref{eq:3.2} imply that
\begin{equation}
  \label{eq:3.3a}
  \begin{split}
& \|\dot{\theta}-E\|_{L^1(0,T)}+\|\dot{E}\|_{L^1(0,T)}
\\ \le &
C(\bM)(\|\phi_{1,E(t)}v^2\|_{L^1(0,T;L^1_x)}
+\|\phi_{2,E(t)}v^2\|_{L^1(0,T;L^1_x)})
\\ & +
C(\bM)(\|\phi_{1,E(t)}f(v)\|_{L^1(0,T;L^1_x)}
+\|\phi_{2,E(t)}f(v)\|_{L^1(0,T;L^1_x)})
\\ \le &
C(\bM)\left(\sum_{i=1,2}\left\|\la x\ra^{2s}\phi_{i,E(t)}
\right\|_{L^\infty(0,T;L^\infty_x)}\right)
\|v\|_{L_t^2(0,T;H_x^{1,-s})}^2
\\ \le & C(\bM)(\bM_2(T)+\bM_3(T))^2,
  \end{split}
\end{equation}
and
\begin{equation}
  \label{eq:3.3b}
  \begin{split}
\|\dot{\theta}-E\|_{L^\infty(0,T)}+\|\dot{E}\|_{L^\infty(0,T)}
\lesssim & \sup_{0\le t\le T}(\|v\|_{H^1}^2+\|v\|_{H^1}^p)
\\ \le & C(\bM)(\bM_4(T)+\bM_5(T))^2.
  \end{split}
\end{equation}
Hereafter we denote by $C(\bM)$ various functions of
$\bM_1$, \dots, $\bM_5$ that are bounded in a finite neighborhood of $0$.
By \eqref{eq:3.3ap} and \eqref{eq:3.3a},
\begin{equation}
  \label{eq:3.4}
\bM_1(T) \lesssim \|u_0\|_{H^1}+ C(\bM)(\bM_2+\bM_3)^2.
\end{equation}

From Remark \ref{rem:1} and \eqref{eq:2.1}, it follows that
\begin{align*}
  |\la w(t), \phi_*\ra|\le & \|v\|_{L^{2,-s}_x}
\sum_{i=1,2}\|\la x\ra^{s}(\phi_{i,E}-\phi_*)\|_{L^2}
\\ \lesssim & |E(t)-E_*|\|w\|_{L_x^{2,-s}},
\end{align*}
and that
\begin{equation}
  \label{eq:3.6}
  \begin{split}
\bM_3(T) \le C(\bM)\bM_1(T)(\bM_2(T)+\bM_3(T)).
  \end{split}
  \end{equation}
Similarly, we have
\begin{equation}
    \label{eq:3.7}
\bM_5(T) \le C(\bM)\bM_1(T)(\bM_4(T)+\bM_5(T)).      
\end{equation}

Next, we will estimate $\bM_2(T)$.
By \eqref{eq:v'},
$$\bM_2(T)\le I_1+I_2+I_3+I_4,$$
where
\begin{align*}
 & I_1=\|e^{-it L}P_cw(0)\|_{L_t^2(0,T;H_x^{1,-s})}, \\
 & I_2=\left\|\int_0^t e^{-i(t-s)L}P_cg_2(s)ds
\right\|_{L^2_t(0,T;H_x^{1,-s})},\\
 & I_3=\left\|\int_0^t e^{-i(t-s)L}P_cf(v(s))ds
\right\|_{L^2_t(0,T;H_x^{1,-s})},\\
 & I_4=\left\|\int_0^t e^{-i(t-s)L}P_c\tilde{g}(s)ds
\right\|_{L^2_t(0,T;H_x^{1,-s})},
\end{align*}
and $\tilde{g}(s)=g_3(s)+g_4(s)-f(v(s))$.
Lemma \ref{lem:2.2} yields $$I_1\lesssim  \|w(0)\|_{H^1}.$$
 By Lemma \ref{lem:2.3}, \eqref{eq:3.3a} and \eqref{eq:3.3b},
\begin{align*}
 I_2 & \lesssim \|P_cg_2\|_{L_t^2(0,T;H_x^{1,s})}
\\ \le &
\left\|P_c\phi_{E(t)}\right\|_{L^\infty(0,T;H_x^{1,s})}
\|\dot{\theta}-E\|_{L^2(0,T)}
+\left\|P_c\pd_E\phi_{E(t)}\right\|_{L^\infty(0,T;H_x^{1,s})}
\|\dot{E}\|_{L^2(0,T)}
\\ \le &
C(\bM)\bM_1(T)^{1/(p-1)}(\bM_2(T)+\bM_3(T)+\bM_4(T)+\bM_5(T))^2.
\end{align*}
Note that $\|P_c\pd_E\phi_E\|_{H^1}\lesssim |E-E_*|^{1/(p-1)}$
follows from Proposition \ref{prop:1}.
By Minkowski's inequality and Lemma \ref{lem:2.2},
\begin{align*}
  I_3 \lesssim & 
\int_0^T\|e^{-i(t-\tau)L}P_cf(v(\tau))\|_{L^2_t(0,T;H_x^{1,-s})}d\tau
\\ \lesssim & \int_0^T\|f(v(s))\|_{H_x^1}ds
\\ \lesssim &
\|v\|^q_{L^q(0,T;W_x^{1,2p})}\|v\|_{L_t^\infty(0,T;H_x^1)}^{p-q},
\end{align*}
where $2/q+1/p=1$. Note that $p\ge q>2$ if $p\ge 3$.
Thus we have
$$I_3\le C(\bM)(\bM_4(T)+\bM_5(T))^p.$$
Since $\tilde{g}=O(\phi_E^{p-1}|v|+\phi_E|v|^{p-1})$,
Lemma \ref{lem:2.3} yields that
\begin{equation}
  \label{eq:I4}
\begin{split}
I_4\lesssim & \|\tilde{g}\|_{L^2_t(0,T;H_x^{1,s})}
\\ \lesssim &
\|\la x\ra^{2s}\phi_{E(t)}^{p-1}\|_{L_t^\infty(0,T;W_x^{1,\infty})}
\|v\|_{L_t^2(0,T;H_x^{1,-s})}
\\ & +\|\la x\ra^s \phi_{E(t)}\|_{L_x^\infty(0,T;W_x^{1,\infty})}
\||v|^{p-1}\|_{L^2_t(0,T;H_x^1)}.
\end{split}  
\end{equation}
Since
$$\||v|^{p-1}\|_{H^1}\le
\|v\|_{W^{1,2(p-1)/(p-2)}}\||v|^{p-2}\|_{L^{2(p-1)}}
\lesssim \|v\|_{W^{1,2(p-1)/(p-2)}}^{p-1},
$$
it follows from \eqref{eq:I4}, Proposition \ref{prop:1}
 and the interpolation theorem that
\begin{align*}
I_4 \le & C(\bM)\left(\|v\|_{L_t^2(0,T;H_x^{1,-s})}+
\|v\|_{L_t^{2(p-1)}\left(0,T;W_x^{1,2(p-1)/(p-2)}\right)}^{p-1}\right)
\\ \le & C(\bM)\{\bM_1(T)(\bM_2(T)+\bM_3(T))
+(\bM_4(T)+\bM_5(T))^{p-1}\}. 
\end{align*}
Combining the above, we see that
\begin{equation}
  \label{eq:3.8}
  \bM_2(T)\le C(\bM)\sum_{1\le i\le 5}\bM_i(T)^2.
\end{equation}

Finally, we will estimate $\bM_4(T)$.
In view of \eqref{eq:v'},
$$\bM_4(T)\le J_1+J_2+J_3,$$
where
\begin{align*}
J_1=&\left\| e^{-itL}P_cw(0)
\right\|_{L^\infty(0,T;H^1_x)\cap L^q(0,T;W_x^{1,2p})}
\\  J_2=&\left\|\int_0^t e^{-i(t-s)L}P_cg_2(s)ds
\right\|_{L^\infty(0,T;H^1_x)\cap L^q(0,T;W_x^{1,2p})},
\\ J_2=&\left\|\int_0^t e^{-i(t-s)L}P_c(g_3(s)+g_4(s))ds
\right\|_{L^\infty(0,T;H^1_x)\cap L^q(0,T;W_x^{1,2p})}.
\end{align*}
Using the Strichartz estimate (Lemma \ref{lem:2.1}), we have
\begin{align*}
J_1 \lesssim & \|w(0)\|_{H^1},\\
J_2 \lesssim & \|P_cg_2(s)\|_{L_t^1(0,T;H^1_x)}ds
\\ \lesssim &
\|\dot{\theta}-E\|_{L^1(0,T)}\sup_{t\in[0,T]}\|P_c\phi_{E(t)}\|_{H^1_x}
+\|\dot{E}\|_{L^1(0,T)}\sup_{t\in[0,T]}\|P_c\pd_E\phi_{E(t)}\|_{H^1_x}.
\end{align*}
Hence by \eqref{eq:3.3a},
$$J_2 \le C(\bM)(\bM_2(T)^2+\bM_3(T)^2).$$
Using the Strichartz estimate  and Corollary
\ref{cor:2.4}, we have
$$
J_3\lesssim \|g_3+g_4\|_{L_t^1(0,T;H_x^1)+L_t^2(0,T;H_x^{1,s})}.$$
Since $g_3(t)+g_4(t)=O(\phi_{E(t)}^{p-1}|v|+|v|^p)$,
\begin{align*}
& \|g_3+g_4\|_{L_t^1(0,T;H_x^1)+L_t^2(0,T;H_x^{1,s})}
\\ \lesssim &
\|\phi_{E(t)}^{p-1}v\|_{L_t^2(0,T;H_x^{1,s})}+\|f(v)\|_{L_t^1(0,T;H_x^1)}
\\ \lesssim &
\|\la x\ra^{2s}\phi_{E(t)}^{p-1}\|_{L_t^\infty(0,T;W_x^{1,\infty})}
\|v\|_{L_t^2(0,T;H_x^{1,-s})}
+\|v\|_{L_t^q(0,T;W_x^{1,2p})}^{q}\|v\|_{L_t^\infty(0,T;H_x^1)}^{p-q}.
\end{align*}
Thus we have
$$
J_3\le C(\bM)\{\bM_1(T)(\bM_4(T)+\bM_5(T))+(\bM_4(T)+\bM_5(T))^p\}.
$$
Combining the above, we have
\begin{equation}
  \label{eq:3.10}
  \bM_4(T)\le C(\bM)\sum_{1\le i\le 5}\bM_i(T)^2.
\end{equation}
\par

It follows from \eqref{eq:3.4}--\eqref{eq:3.7}, \eqref{eq:3.8} and
\eqref{eq:3.10} that if $\varepsilon_0$ is sufficiently small,
\begin{equation}
  \label{eq:3.12}
\sum_{1\le i\le 5}\bM_i(T)\lesssim \|u_0\|_{H^1}.  
\end{equation}
Thus by continuation argument, we may let $T\to\infty$.
By \eqref{eq:3.3a}, there exists an $E_+<0$ satisfying
$\lim_{t\to\infty}E(t)=E_+$ and $|E_+-E_*|\lesssim \|u_0\|_{H^1}.$
In view of \eqref{eq:3.12}, we see that
$$
w_1:=-i\lim_{t\to\infty}\sum_{2\le j\le 4}\int_0^t e^{isL}P_c
e^{-i\theta(s)}g_j(s)ds$$
exists in $H^1$ and that
\begin{align*}
\|w_1\|_{H^1}\lesssim &  \|g_2(s)\|_{L_t^1H_x^1}
+\|g_3+g_4\|_{L_t^2H_x^{1,s}+ L_t^1H_x^1}
\\ \lesssim & \|u_0\|_{H^1},\\
\lim_{t\to\infty}\|P_cw(t)-& e^{-itL}(P_cw(0)+w_1)\|_{H^1}=0.
\end{align*}
By \cite{Sch}, we have
$\|e^{-itL}P_cf\|_{L^4}\lesssim t^{-1/2}\|f\|_{L^{4/3}}.$
Since $L^{4/3}(\R^2)\cap H^1(\R^2)$ is dense in $H^1(\R^2)$,
it follows that
$\|e^{-itL}(P_cw(0)+w_1)\|_{L^{4}}\to0$ as $t\to\infty,$ and that
\begin{equation}
  \label{eq:nyan1}
\begin{split}
& \|P_cw(t)\|_{L^4} \\ \le
&\|P_cw(t)-e^{-itL}(P_cw(0)+w_1)\|_{H^1}
+\|P_ce^{-itL}(P_cw(0)+w_1)\|_{L^4}
\\ \to & 0 \quad\text{as $t\to\infty$.}
\end{split}  
\end{equation}
Analogously to \eqref{eq:3.6}, we have
\begin{equation}
  \label{eq:nyan2}
\|P_{d}w(t)\|_{H^1}\lesssim \|P_{d}w(t)\|_{L^4}
\lesssim |E(t)-E_*|\|P_cw(t)\|_{L^4}.
\end{equation}
Combining \eqref{eq:nyan1} and \eqref{eq:nyan2}, we have
$\lim_{t\to\infty}\|P_{d}w(t)\|_{H^1}=0.$
Thus by \eqref{eq:1.3} and \eqref{eq:trans}, 
$$
\lim_{t\to\infty}\left\|u(t)-e^{-i\theta(t)}\phi_{E(t)}
-e^{-itL}P_c(w(0)+w_1)\right\|_{H^1}=0.$$
Thus we complete the proof of Theorem \ref{thm1}.
\end{proof}

\section{Dispersive estimates}
\label{sec:4}
Let $R(\lambda)=(L-\lambda)^{-1}$ and $dE_{ac}(\lambda)$ be the absolute
continuous part of the spectrum measure. 
By the spectral decomposition theorem, we have
\begin{equation}
  \label{eq:invLap}
\begin{split}
P_ce^{-itL}f=& \int_{-\infty}^\infty e^{-it\lambda}dE_{ac}(\lambda)f
\\ =& \frac{1}{2\pi i}\int_{-\infty}^\infty e^{-it\lambda}
P_c(R(\lambda+i0)-R(\lambda-i0))fd\lambda.
\end{split}  
\end{equation}
We will prove Lemma \ref{lem:2.2} by using Plancherel's theorem and
the following estimate on the resolvent $R(\lambda)$. 
\begin{lemma}
  \label{lem:4.1} Let $s>1$. Then there exists a positive constant $C$
such that
$$\|R(\lambda\pm i0)P_cf\|_{L^2_\lambda(0,\infty;L^{2,-s}_x)}
\le C\|f\|_{L^2}$$
for every $f\in L^2(\R^2)$.
\end{lemma}
First, we prove Lemma \ref{lem:2.2} assuming Lemma \ref{lem:4.1}.
\begin{proof}[Proof of Lemma \ref{lem:2.2}]
  By the inversion of the Laplace formula (see \cite{Pa}), we have
  \begin{align*}
e^{-itL}P_cf=& \frac{1}{2\pi i}\int_{-\infty}^\infty d\lambda e^{-it\lambda}
(R(\lambda+i0)-R(\lambda-i0))P_cf
\\ =& 
\frac{(it)^{-j}}{2\pi i}\int_{-\infty}^\infty d\lambda e^{-it\lambda}
\pd_\lambda^j(R(\lambda+i0)-R(\lambda-i0))P_cf
 \quad\text{in $\mathcal{S}'_x(\R^2)$} 
  \end{align*}
for any $t\ne0$ and $f\in \mathcal{S}_x(\R^2)$.  
Since
$$\|\pd_\lambda^jR(\lambda\pm i0)P_c
\|_{B(L^{2,(j+1)/2+0},L^{2,-(j+1)/2-0})}
\lesssim \la \lambda\ra^{-(j+1)/2},$$
the above integral absolutely converges in $L^{2,-(j+1)/2-0}_x$ for $j\ge2$. 
\par

Suppose $g(t,x)=g_1(t)g_2(x)$,
$g_1\in C_0^\infty(\R\setminus\{0\})$ and $g_2\in\mathcal{S}(\R^2)$.
Making use of Fubini's theorem and integration by parts,
we have for $j\ge2$,
\begin{align*}
& \la e^{-itL}P_cf,g\ra_{t,x}
\\=&
\frac1{2\pi i}\int_{-\infty}^\infty dt (it)^{-j}g_1(t)
\int_{-\infty}^\infty d\lambda e^{-it\lambda}
\pd_\lambda^j\left\la(R(\lambda+i0)-R(\lambda-i0))P_cf,g_2\right\ra_x
\\=&
\frac{1}{2\pi i}\int_{-\infty}^\infty d\lambda
\pd_\lambda^j \left\la(R(\lambda+i0)-R(\lambda-i0))P_cf,g_2 \right\ra_x
\int_{-\infty}^\infty dt (it)^{-j}g_1(t) e^{-it\lambda}
\\=& \frac{1}{\sqrt{2\pi}i}\int_{-\infty}^\infty d\lambda
(\mathcal{F}_tg_1)(\lambda)
\left\la(R(\lambda+i0)-R(\lambda-i0))P_cf,g_2\right\ra_x.
\end{align*}
Hence it follows from the above that
\begin{equation*}
  \left\la e^{-itL}P_cf,g\right\ra=\frac{1}{\sqrt{2\pi}i}
\int_{-\infty}^\infty d\lambda \left\la (R(\lambda+i0)-R(\lambda-i0))P_cf,
\mathcal{F}_t g(\lambda,\cdot) \right\ra_x
\end{equation*}
for every $g\in C_0^\infty(\R_t\setminus\{0\})\otimes\mathcal{S}(\R^2_x)$.
Using Plancherel's theorem, we have
\begin{equation}
  \label{eq:pl1}
  \begin{split}
& \left|\la e^{-itL}P_cf,g\ra_{t,x}\right|
\\ \le & \frac{1}{\sqrt{2\pi}}\int_{-\infty}^\infty d\lambda
\|(R(\lambda+i0)-R(\lambda-i0))P_cf\|_{L^{2,-s}_x}
\|\mathcal{F}_t g(\lambda,\cdot) \|_{L^{2,s}_x}
\\ \le &
(2\pi)^{-1/2}
\|(R(\lambda+i0)-R(\lambda-i0))P_cf\|_{L^2_\lambda(0,\infty;L^{2,-s}_x)}
\|g\|_{L_t^2L_x^{2,s}}.    
  \end{split}
\end{equation}
Since $C_0^\infty(\R_t\setminus\{0\})\otimes\mathcal{S}(\R^2_x)$
is dense in $L^2_tL_x^{2,s}$, Lemma \ref{lem:2.2} immediately follows from
\eqref{eq:pl1}.
\end{proof}
\par

Now, we turn to prove Lemma \ref{lem:4.1}. 
First, we will investigate the free resolvent operator
$R_0(\lambda)$ in $\R^2$.
\begin{lemma}
  \label{lem:4.2}
There exists a positive constant $C$ such that
$$\sup_x\|R_0(\lambda\pm i0)f\|_{L^2_\lambda(0,\infty)}
\le C\|f\|_{L^2}$$ for every $f\in L^2(\R^2)$.
\end{lemma}
\begin{remark}
Obviously, the estimate
$\|R_0(\lambda\pm i0)\|_{B(L^{2,s,},L^{2,-s})}
\lesssim\la\lambda\ra^{-1/2}$ does not suffice to prove
Lemma \ref{lem:4.2}. We will use the
boundedness of the Hankel transform in $L^2_{rad}$.  
\end{remark}

\begin{proof}[Proof of Lemma \ref{lem:4.2}]
For any  $k\ge 0$, we have
$$R_0(k^2\pm i0)f(x)=\frac{\pm i}{4}\int_{\R^2}H_0(k|x-y|)f(y)dy,$$
where $H_0^\pm$ are the Hankel functions of order $0$ and
$$H_0^\pm(z)=J_0(z)\pm Y_0(z).$$

Let $(\tau_xf)(y):=f(x-y)$ and decompose $\tau_xf\in L^2(\R^2)$ 
into a Fourier series as
$$\tau_xf=\sum_{m\in\Z}f_{x,m}(r)e^{im\theta}\in \bigoplus_{m\in \Z}
e^{im\theta}L_{rad}.$$
Then
\begin{align*}
R_0(k^2\pm i0)f(x)
=& \frac{\pm i}{4}\int_{\R^2}H_0^{\pm}(k|y|)\tau_xf(y)dy
\\=& 
\frac{\pm\pi i}{2}\int_0^\infty H_0^\pm(kr)f_{x,0}(r)rdr.
\end{align*}
Titchmarsh \cite{Tit} and Rooney \cite{R3} tell us that
the operators $T_1$ and $T_2$ defined by
$$T_1f(x)=\int_0^\infty J_0(xy)f(y)dy,\quad
T_2f(x)=\int_0^\infty Y_0(xy)f(y)dy,$$
are bounded on $L^2_{rad}$. Thus we have
$$\sup_x\left(\int_0^\infty |R_0(k^2\pm i0)f|^2kdk\right)^{1/2}
\lesssim \|f_{x,0}\|_{L^2_{rad}}^2.$$
Since
$$\|f\|_{L^2}=\|\tau_xf\|_{L^2}=\left(2\pi\sum_{m\in \Z}
\int_0^\infty |f_{x,m}(r)|^2rdr\right)^{1/2},$$
it follows that
$$\sup_x\|R_0(\lambda\pm i0)f\|_{L^2_\lambda(0,\infty)}
\lesssim \|f\|_{L^2}.$$ 
Thus we complete the proof of Lemma \ref{lem:4.2}.
\end{proof}

We will prove Lemma \ref{lem:4.1} by using Lemma \ref{lem:4.2} and
the resolvent expansion obtained by Schlag \cite{Sch} based on
Jensen and Nenciu \cite{JeN}.
\par

Before we prove Lemma \ref{lem:4.1}, let us introduce a definition of
the non-resonance condition given by Jensen and Nenciu \cite{JeN}.
\begin{definition}\label{def:NR}
Let $v(x)=|V(x)|^{1/2}$ and let $P$ and $Q$ be orthogonal projections
defined by
$$Pf=\frac{\la f,v\ra v}{\|V\|_{L^1}}, \quad Q=I-P.$$
We say that $0$ is not a resonance of $L$ if   
$D_0:=Q(U+vG_0v)Q$ is invertible on $QL^2(\R^2)$.
\end{definition}

\begin{proof}[Proof of Lemma \ref{lem:4.1}]
For every $f\in\mathcal{S}(\R^2)$, we have
\begin{equation}
  \label{eq:4.1}
  R(\lambda\pm i0)f=R_0(\lambda\pm i0)f-R_0(\lambda\pm i0)VR(\lambda\pm i0)f.
\end{equation}
By Lemma \ref{lem:4.2}, there exists a $C>0$ such that
for every $f\in L^2(\R^2)$,
\begin{equation}
  \label{eq:4.1a}
\begin{split}
\|R_0(\lambda\pm i0)f\|_{L_x^{2,-s}L_\lambda^2(0,\infty)}  \le &
\|\la x\ra^{-s}\|_{L^2}
\|R_0(\lambda\pm i0)f\|_{L_x^{\infty}L_\lambda^2(0,\infty)}
\\ \le & C\|f\|_{L^2}.
\end{split}  
\end{equation}
\par

Next, we deal with the low energy part of the second term of
\eqref{eq:4.1}. As \cite{JeN,Sch}, we put $U(x)=1$ for
$x\in V^{-1}([0,\infty))$, $U(x)=-1$ for $x\in V^{-1}((-\infty,0))$, and 
$M^\pm(\lambda):=U+vR_0(\lambda\pm i0)v$. Then
$$ R_0(\lambda\pm i0)VR(\lambda\pm i0)f
=R_0(\lambda\pm i0)vM^\pm(\lambda)^{-1}vR_0(\lambda\pm i0)f.$$
Schlag \cite[Lemma 9]{Sch} tells us that
\begin{equation}
  \label{eq:4.2}
  M^\pm(\lambda)^{-1}=h_\pm(\lambda)^{-1}S+QD_0Q+E^\pm(\lambda)
\quad\text{in $B(L^2(\R^2))$,}
\end{equation}
where $S$ is a finite rank operator,
$\|E^\pm(\lambda)\|_{B(L^2)}=O(\lambda^{1/4})$ as $\lambda\to0$, and
\begin{equation}
  \label{eq:4.3}
h_+(\lambda)=a\log\lambda+z,\quad h_-(\lambda)=\overline{h_+(\lambda)},
\end{equation}
and $a\in\R$ and $z\in \C$ are constants with $a\ne0$ and $\Im z\ne0$.
\par

Let $\lambda_1$ be a sufficiently small positive number. From
\cite[Lemma 5]{Sch}, it follows that for $0<\lambda\le\lambda_1$,
\begin{equation}
  \label{eq:4.4}
R_0(\lambda\pm i0)=c_\pm(\lambda)P_0+G_0+E_0^\pm(\lambda)
\quad\text{in $B(L^{2,s},L^{2,-s})$,}
\end{equation}
and
\begin{equation}
\label{eq:4.4b}
  \|E_0^\pm(\lambda)\|_{B(L^{2,s},L^{2,-s})}=O(\lambda^{1/4}),
\end{equation}
where $P_0f=\la f,1\ra_x$, $G_0=(-\Delta)^{-1}$, $\gamma$ is the Euler
number and 
\begin{equation}
  \label{eq:4.5}
c_\pm(\lambda)=\pm\frac{i}4-\frac{\gamma}{2\pi}-\frac{1}{4\pi}\log
\left(\frac{\lambda}{4}\right).
\end{equation}
Thus $\widetilde{R}_0^\pm(\lambda)=R_0(\lambda\pm i0)-c_\pm(\lambda)P_0$
satisfies
\begin{equation}
  \label{eq:4.6}
\sup_{0<\lambda<\lambda_1}\|\widetilde{R}_0^\pm(\lambda)
\|_{B(L^{2,s},L^{2,-s})} <\infty.
\end{equation}
\par

Let $\chi(\lambda)$ be a characteristic function on $[0,\lambda_1]$.
 Using Lemma \ref{lem:4.2}, \eqref{eq:4.2}, \eqref{eq:4.6}
and the fact that
$v(x)\lesssim \la x\ra^{-\sigma/2}$ with $\sigma>3$,
we have
 \begin{align*}
& \| \chi(\lambda) \widetilde{R}_0^\pm(\lambda)
vM^\pm(\lambda)^{-1}v R_0(\lambda\pm i0)f
\|_{L_\lambda^2(0,\infty;L_x^{2,-s})}
\\ \le &   \sup_{0<\lambda<\lambda_1}
\|\widetilde{R}_0^\pm(\lambda)\|_{B(L^{2,s},L^{2,-s})}
\left\|\|\chi(\lambda)vM^\pm(\lambda)^{-1}v R_0(\lambda\pm i0)f
\|_{L_x^{2,s}} \right\|_{L^2_\lambda(0,\infty)}
\\ \lesssim & \|\chi(\lambda)vR_0(\lambda\pm i0)f\|_{L^2_{x,\lambda}}
\\ \lesssim & 
\|v\|_{L^2_x}\sup_x\|R_0(\lambda\pm i0)f\|_{L_\lambda^2(0,\infty)}
\lesssim \|f\|_{L^2}
 \end{align*}
for any $s\in(1,3/2)$. 
Since $P_0vQ=0$, it follows from \eqref{eq:4.2} that
$$
 c^\pm(\lambda)P_0vM^\pm(\lambda)^{-1}vR_0(\lambda\pm i0)
=I_1+I_2,$$
where
\begin{align*}
I_1=& c^\pm(\lambda)h_\pm(\lambda)^{-1}P_0vSvR_0(\lambda\pm i0),\\ 
I_2=& c^\pm(\lambda)h_\pm(\lambda)^{-1}P_0vE^\pm(\lambda)
vR_0(\lambda\pm i0).
\end{align*}
By \eqref{eq:4.3}, \eqref{eq:4.5}, 
$\sup_{0<\lambda\le\lambda_1}|c_\pm(\lambda)/h_\pm(\lambda)|<\infty.$
Hence  from Lemma \ref{lem:4.1},
\begin{align*}
  \|I_1f\|_{L_\lambda^2(0,\infty;L_x^{2,-s})}\le &
\|\la x\ra^{-s}\|_{L^2}
\|vSvR_0(\lambda\pm i0)f\|_{L_\lambda^2(0,\infty;L_x^1)}
\\ \lesssim & \|v\|_{L^2}
\|vR_0(\lambda\pm i0)f \|_{L_\lambda^2(0,\infty;L_x^2)}
\\ \lesssim &  \|v\|_{L^2}^2
\sup_x\|R_0(\lambda\pm i0)f\|_{L_\lambda^2(0,\infty)}
\\ \lesssim & \|f\|_{L^2}.
\end{align*}
Using Schwarz's inequality and \eqref{eq:4.4b}, we have
\begin{align*}
  \|P_0vE^\pm(\lambda)vR_0(\lambda\pm i0)f\|_{L_x^{2,-s}}
\lesssim & \|v\|_{L^2}\|E^\pm(\lambda)vR_0(\lambda\pm i0)f\|_{L^2}
\\ \lesssim & |\lambda|^{1/4}\|vR_0(\lambda\pm i0)f\|_{L^2}.
\end{align*}
Hence it follows that
\begin{align*}
  \|\chi(\lambda)I_2\|_{L^2(0,\infty;L_x^{2,-s})} \lesssim &
\sup_{\lambda>0}\left(\chi(\lambda)|\lambda|^{1/4}|c_\pm(\lambda)|
\right) \|vR_0(\lambda\pm i0)f\|_{L_{x,\lambda}^2}
\\ \lesssim & \sup_x\|R_0(\lambda\pm i0)f\|_{L_\lambda^2}
\\ \lesssim & \|f\|_{L^2}.
\end{align*}
Combining the above, we obtain
\begin{equation}
\label{eq:low}
\|\chi(\lambda)R_0(\lambda\pm i0)VR(\lambda\pm i0)f
\|_{L_\lambda^2(0,\infty;L_x^{2,-s})} \lesssim \|f\|_{L^2}.
\end{equation}
\par

Next, we consider the high energy part.
The assumptions (H2) and (H3) imply that
\begin{equation}
  \label{eq:4.9}
\sup_{\lambda\ge \lambda_1}\|R(\lambda\pm i0)P_c\|_{B(L^{2,s},L^{2,-s})}
\lesssim \la \lambda\ra^{-1/2},  
\end{equation}
See \cite[Appendix A]{Ag} and \cite{Mu} for the proof.
Let $\tilde{\chi}(\lambda)=1-\chi(\lambda)$. By \eqref{eq:4.9} and 
Fubini's theorem,
\begin{equation}
  \label{eq:4.10}
\begin{split}
& \|\tilde{\chi}(\lambda)P_cR(\lambda\pm i0)VR_0(\lambda\pm i0)f
\|_{L_\lambda^2(0,\infty;L_x^{2,-s})}
\\ \lesssim  & \bigl\|\|VR_0(\lambda\pm i0)f
\|_{L_x^{2,s}}\bigr\|_{L_\lambda^2(0,\infty)}
\\ \le & \|V\|_{L^{2,s}}
\sup_x\|R_0(\lambda\pm i0)f\|_{L_\lambda^2(0,\infty)}
\lesssim  \|f\|_{L^2}.
\end{split}
\end{equation}
Combining \eqref{eq:4.1}, \eqref{eq:4.1a}, \eqref{eq:low} and
\eqref{eq:4.10}, we obtain
$$\|R(\lambda\pm i0)P_cf\|_{L^2_\lambda(0,\infty;L_x^{2,-s})}
\le C\|f\|_{L^2}.$$
Thus we complete the proof of Lemma \ref{lem:4.1}.
\end{proof}
\par
Next, we will prove Lemma \ref{lem:2.3}.
For the purpose, we need the following.
\begin{lemma}
  \label{lem:4.3} Assume \textrm{(H1)}--\textrm{(H3)}.
Let $g(t,x)\in \mathcal{S}_t(\R)\otimes\mathcal{S}_x(\R^2)$
and
$$U(t,x)=\frac{1}{\sqrt{2\pi}i}\int_{-\infty}^\infty d\lambda
e^{-it\lambda}
\{R(\lambda-i0)+R(\lambda+i0)\}P_c(\mathcal{F}^{-1}_tg)(\lambda,\cdot).
$$
Then,
\begin{align*}
U(t,x)=& 2\int_0^t ds e^{-i(t-s)L}P_cg(s,\cdot)
+\int_{-\infty}^0 ds e^{-i(t-s)L}P_cg(s,\cdot)
\\ & -\int_0^\infty ds e^{-i(t-s)L}P_cg(s,\cdot).
\end{align*}
\end{lemma}
\begin{proof}
Since Lemma \ref{lem:4.3} can be proved in the same as that of
Lemma 11 in \cite{Mi}, we omit the proof.
\end{proof}
\begin{proof}[Proof of Lemma \ref{lem:2.3}]
Suppose that $g(t,x)$ and $h(t,x)$ belong to
$\mathcal{S}_t(\R)\otimes\mathcal{S}_x(\R^2)$. 
It follows from Fubini's theorem that
\begin{align*}
&  \la U, h\ra_{t,x}\\=& \frac{1}{\sqrt{2\pi}i}
\int_{-\infty}^\infty d\lambda \int_{-\infty}^\infty dt e^{-it\lambda}
\left\la (R(\lambda+i0)+R(\lambda-i0))P_c 
\mathcal{F}_t^{-1}g(\lambda,\cdot),h(t,\cdot)\right\ra_x
\\=&i^{-1}
\int_{-\infty}^\infty d\lambda
\left\la (R(\lambda+i0)+R(\lambda-i0))P_c
\mathcal{F}_t^{-1}g(\lambda,\cdot),
\mathcal{F}_th(\lambda,\cdot)\right\ra_x.
\end{align*}
Using Plancherel's theorem and \eqref{eq:4.9}, we obtain
\begin{align*}
& \left|\la U,h\ra_{t,x}\right| \\ \le &
\left\|(R(\lambda+i0)+R(\lambda-i0))P_c
\mathcal{F}_t^{-1}g(\lambda,\cdot)\right\|_{L_\lambda^2L_x^{2,-s}}
\left\|\mathcal{F}_th(\lambda,\cdot)\right\|_{L_\lambda^2L_x^{2,s}}
\\ \le &
\sup_{\lambda\in\R}
\left\|(R(\lambda+i0)+R(\lambda-i0))P_c\right\|_{B(L^{2,s},L^{2,-s})}
\|g\|_{L_t^2L_x^{2,-s}}\|h\|_{L_t^2L_x^{2,s}}
\end{align*}

Since $\mathcal{S}_t(\R)\otimes\mathcal{S}_x(\R^2)$ is dense
in $L^2_tL_x^{2,s}$ and $L_t^2L_x^{2,-s}$, we see that
\begin{equation}
 \label{eq:U-est}
\|U\|_{L_t^2L_x^{2,-s}}\lesssim \|g\|_{L_t^2L_x^{2,s}}
\end{equation}
holds for every $g\in L^2_tL_x^{2,s}$.
\par
On ther other hand, Lemma \ref{lem:2.2} implies
\begin{align*}
 \left\|\int_I e^{-i(t-s)L} Qg(s)ds\right\|_{L_t^2L_x^{2,-s}}
 \lesssim \left\|\int_{I} e^{isL}g(s)ds\right\|_{L^2}
 \lesssim  \|g\|_{L_t^2L_x^{2,s}}
\end{align*}
for every $g\in L_t^2L_x^{2,s}$ and $I\subset \R$.
 Combining the above with \eqref{eq:U-est} and Lemma \ref{lem:4.3},
we obtain Lemma \ref{lem:2.3}. Thus we complete the proof.
\end{proof}
\par

Finally, we prove Corollary \ref{cor:2.4}.
\begin{proof}[Proof of Corollary \ref{cor:2.4}]
Let $(q,r)$ be admissible and let $T$ be an operator defined by
$$Tg(t)=\int_\R ds e^{-i(t-s)L}P_cg(s).$$
Lemmas \ref{lem:2.1} and \ref{lem:2.2} yield
$f:=\int_\R e^{isL}P_cg(s)ds \in L^2(\R)$  and that there exists a $C>0$
such that
\begin{equation}
  \label{eq:T}
\|Tg(t)\|_{L^q_tL^r_x}\le C \|g\|_{L_t^2L_x^{2,s}}
\end{equation}
for every $g\in L_t^2L_x^{2,s}$.
Since $q>2$, it follows from Lemma 3.1 in \cite{SmS} and \eqref{eq:T}
that
\begin{equation}
  \label{eq:T1}
\left\|\int_{s<t} ds e^{-i(t-s)L}P_cg(s)\right\|_{L_t^qL_x^r} \lesssim
\|g\|_{L_t^2L_x^{2,s}}.
\end{equation}
Thus we prove Corollary \ref{cor:2.4}.
\end{proof}

{\scshape
\begin{flushright}
\begin{tabular}{l}
Faculty of Mathematics \\
Kyushu University \\
6-10-1 Hakozaki\\
Fukuoka 812-8581, Japan\\
{\upshape e-mail: mizumati@math.kyushu-u.ac.jp}\\
\end{tabular}
\end{flushright}
}
\end{document}